\documentclass{article}
\usepackage{amsmath}
\usepackage{amssymb}

\newcommand{\EE}{\mathbb{E}}

\newcommand{\NN}{\mathbb{N}}
\newcommand{\PP}{\mathbb{P}}
\newcommand{\RR}{\mathbb{R}}
\newcommand{\SaS}{\mathbb{S}}
\newcommand{\T}{\mathbb{T}}
\newcommand{\bP}{{\bf P}}

\newcommand{\Aa}{{\cal A}}
\newcommand{\C}{{\cal C}}
\newcommand{\hH}{{\cal H}}
\newcommand{\I}{{\cal I}}

\newcommand{\U}{{\cal U}}

\newcommand{\Int}{{\rm int}}

\newtheorem{theorem}{Theorem}
\newtheorem{lemma}{Lemma}

\newtheorem{remark}{Remark}

\begin{document}
\begin{center}
{\Large The $\beta$-Mixing Rate of STIT Tessellations}

\vspace{1cm}

{Servet Mart{\'i}nez}\\
{Departamento Ingenier{\'\i}a Matem\'atica and Centro
Modelamiento Matem\'atico,\\  Universidad de Chile,\\
UMI 2807 CNRS, Casilla 170-3, Correo 3, Santiago, Chile.\\
Email: smartine@dim.uchile.cl} 

\vspace{0.5cm}

{Werner Nagel}\\
{Friedrich-Schiller-Universit\"at Jena,\\
Institut f\"ur Stochastik,\\
Ernst-Abbe-Platz 2,
D-07743 Jena, Germany.\\
Email: werner.nagel@uni-jena.de} 

\end{center}

\vspace{1cm}

\begin{abstract}
\noindent We consider homogeneous STIT tessellations $Y$ in the $\ell$-dimensional 
Euclidean space $\RR^\ell$ and show that the (spatial) $\beta$-mixing rate 
converges to zero. 
\end{abstract}


\vspace{1cm}


\section{Introduction}

In stochastic geometry, ergodic and mixing properties as well as
weak dependencies in space are studied.
For models which are $\beta$-mixing (or absolutely regular), 
Laws of Large Numbers and Central Limit 
Theorems can be derived, see \cite{hein94, heinschm, heinmol}. 
Moreover, in \cite{heinpaw} the $\beta$-mixing (absolute regularity) 
rate is used to provide a bound for the total variation distance between 
the (reduced) Palm distribution 
and the distribution of a homogeneous point process.
Recently, \cite{calchen13, chen13} formulated sufficient 
weak dependency conditions which allow to derive results for 
order statistics for certain functionals of cells of tessellations.
All these results refer somehow to an underlying point process 
(the point process itself or a germ-grain-model or the Voronoi-tessellation) 
and make use of the mixing conditions of this point process.

\medskip

Because the STIT tessellations are essentially defined and constructed by 
cell division, using lines (in the plane) or hyperplanes (in the general case) 
respectively, the methods used for the above mentioned point-process-based 
models cannot be applied for STIT. In our recent paper \cite{martnag14} 
we developed the concept of encapsulation, and this allowed us to prove that 
the tail-$\sigma$-algebra is trivial. 
In fact, there was also shown the $\alpha$-mixing property. 
On that base we will show now the stronger condition that 
the $\beta$-mixing rate (in space) converges to zero. 

\medskip

The concepts of ergodicity, some weak mixing (in the ergodic-theoretic sense) and tail triviality for stationary 
random measures (including point processes) were already dealt with in \cite{dvj}.
On the other hand, for sequences of random variables,
there are well established definitions of various mixing rates 
(see e.g. the survey paper \cite{bradsurvey}).
But for random spatial models (like point processes, germ-grain models, 
random tessellations) alternative concepts of mixing had to 
be developed, which form an  appropriate base for proving limit theorems. 
A seminal paper for this is \cite{hein94}. 
There a pair of windows $W'=[-a,a]^\ell \subset W= [-b,b]^\ell$, $0<a<b$, 
in the $\ell$-dimensional Euclidean space $\RR^\ell$ 
is chosen and the degree of dependency of the events inside $W'$ on the 
events outside $W$ is used to define mixing rates. For the corresponding 
$\sigma$-algebras 
(referring to the interior of $W'$ and to the complement $W^c$, respectively) 
the various mixing rates are defined in an analogous way as for sequences 
of random variables.

\medskip

In Section \ref{sectionbeta} of the present paper the definition of the 
$\beta$-mixing rate on the space $\T$ of tessellations 
in the $\ell$-dimensional Euclidean space $\RR^\ell$ is given. 
Further, we adapt to the tessellation model a condition, which is also called 'weak Bernoulli condition' 
which is equivalent to the requirement that 
$\beta (a,b) \to 0$  as $b\to \infty$ (see \cite{bradsurvey}, and 
\cite{bradI}, p. 114, for a proof). 
This condition will later be used to prove our 
main result.

\medskip

In Section \ref{sectionSTIT} we recall the construction of STIT 
tessellations as well as some essential properties. In Section 
\ref{sephyper} we supply the result for sets of hyperplanes separating 
facets and state one of the main requirements on $\Lambda$.

\medskip

Our main results are Theorems \ref{TheorSTITbeta} and \ref{eqivbetamix}.
For the distribution of a STIT tessellation and all $a>0$, the $\beta$-mixing 
rate $\beta (a,b)$ converges to zero for $b\to \infty$. 
A detailed proof is given in Section \ref{sectionproofs}. 
Note that in \cite{hls} a definition of $\beta$-mixing (or absolute regularity) 
is introduced which requires not only a convergence 
to zero but, moreover, certain rates of decay, for $a<b<2a$ and $b>2a$ 
respectively. In the present paper we do not prove $\beta$-mixing in this 
stronger sense, but we provide a formula of an upper bound for $\beta$ 
which probably can be sharpened and exploited further.

\medskip

\section{The $\beta$-mixing rate for random tessellations}\label{sectionbeta}

\subsection{ The measurable space of tessellations}
\label{Sub1.1}
A detailed and sound definition of the measurable space of tessellations of a 
Euclidean space is given in \cite{sw} (Ch. 10, Random Mosaics).
A tessellation is a set $T$ of polytopes (the cells) with 
disjoint interiors and covering the Euclidean space, where each bounded 
subset of  $\RR^\ell$ is intersected by only finitely many cells 
(locally finiteness condition). On the other hand, a tessellation can 
as well be considered as a closed subset $\partial T$ which is 
the union of the cell boundaries. There is an obvious one-to-one 
relation between both ways of description of a tessellation, and their 
measurable structures can be related appropriately, see \cite{sw, martnag}.
Denote by $\T$ the set of all tessellations of $\RR^\ell$.

\medskip

Let ${\cal C}$ be the set of all compact subsets of 
$\RR^\ell$. We endow $\T$ with the Borel $\sigma$-algebra ${\cal B}(\T)$
of the Fell topology (also known as the topology of closed convergence), 
namely 
$$ 
{\cal B}(\T)=\sigma \left( \{ \{ T\in 
\T :\, \partial T\cap C=\emptyset \} :\, C\in {\cal C} \} \right)\, . 
$$
(As usual, for a class of sets $\I$ we denote by $\sigma({\I})$ the smallest 
$\sigma$-algebra containing $\I$.)

\medskip

A compact convex polytope $W$ with non-empty interior in $\RR^\ell$,
is called a {\it window}. We can consider tessellations of $W$ and 
denote the set of all those tessellations by $\T\wedge W$. If
$T\in \T$ we denote by $T\land W$ the induced tessellation on $W$.
Its boundary is defined by 
$\partial (T\wedge W)= (\partial T \cap W)\cup {\partial W}$.

\medskip

For a window $W$ we introduce the notation,
$$ 
{\cal B}(\T_{W})=\sigma \left( \{ \{ T\in \T :\, \partial T\cap C=
\emptyset \} :\, C\subseteq W,\, C\in {\cal C} \} \right) . 
$$ 
By definition ${\cal B}(\T_{W})\subset {\cal B}(\T)$ is a 
sub-$\sigma$-algebra. We notice that if $W'\subseteq W$ then 
${\cal B}(\T_{W'})\subseteq {\cal B}(\T_{W})$.
(We have denoted by $\subset$ strict inclusion,
by $\subseteq$ we mean inclusion.)

\medskip

In order to study the $\beta$-mixing rate, we will also consider sets of 
tessellations which are determined by their behavior outside a window $W$, 
i.e. in its complement $W^c$. We define the $\sigma$-algebra
$$
{\cal B}(\T_{W^c})=\sigma \left( \{ \{ T\in \T :\, \partial T\cap C=
\emptyset  \} :\, C\subset W^c,\, C\in {\cal C}  \} \right) .
$$
We have ${\cal B}(\T_{W^c})\subset {\cal B}(\T )$. On the other hand, if
$W'\subseteq W$ then ${\cal B}(\T_{{W}^c})\subseteq {\cal B}(\T_{W'^c})$.

\subsection{The $\beta$-mixing rate}
\label{Sub1.2}

In \cite{hein94} the $\beta$-mixing rate for a random tessellation with 
distribution $\bP$  on $(\T, {\cal B}(\T))$ is introduced  as follows.
Consider a pair of windows $W'=[-a,a]^\ell \subset W= [-b,b]^\ell$, $0<a<b$. 
Then define
\begin{equation}
\label{defbeta}
\beta (a,b) =  \frac{1}{2}\;  \sup_{({\overline {\cal E}}, {\overline {\Aa}})} 
\;\; \sum_{i=1}^I \sum_{j=1}^J |
\bP( {\cal E}_i \cap \Aa_j)-\bP ( {\cal E}_i) \bP( \Aa_j) | ,
\end{equation}
where the supremum is taken over all pairs of finite partitions of $\T$:
${\overline {\cal E}}= \{{\cal E}_i, i=1,\ldots I \}$ and 
${\overline {\Aa}}=\{ \Aa_j, j=1,\ldots J \}$ with $I,J\in \NN$, 
for events ${\cal E}_i\in {\cal B}(\T_{W^c})$, $\Aa_j\in {\cal B}(\T_{W'})$.

\medskip

The condition (II)(b) in Theorem 3.32 in \cite{bradI} 
(or (2.4) in \cite{bradsurvey})  can be formulated in the 
context of random tessellations as
\begin{eqnarray}
\nonumber
&{}& \forall a>0,\ \forall \epsilon>0,\ 
\exists {\cal D}\in {\cal B}(\T_{W'}) ,\ \bP ( {\cal D})>1-\epsilon ,\ 
\exists b>a, \\ 
&{}& 
\label{bradcond}
\forall {\cal E}\in {\cal B}(\T_{W^c}),\ 
\forall \Aa\in {\cal B}(\T_{W'}) : \\
&{}& 
\nonumber
\Aa\subseteq {\cal D},\ \bP ( \Aa)>0 \Rightarrow 
|\bP( {\cal E}|\Aa)-\bP( {\cal E}) |< \epsilon .
\end{eqnarray} 

Since the proof given in \cite{bradI} refers to families of 
$\sigma$-algebras, it can be applied also to the definitions for random 
tessellations above, and we have immediately the following result:

\begin{theorem}
\label{bradequiv}
The property $\;\lim_{b\to \infty } \beta (a,b) =0$ 
for all $a>0$ is equivalent to condition 
(\ref{bradcond}). \hfill $\Box$
\end{theorem}

\medskip

The following lemma is a corollary of Bradley's 
proof (see \cite{bradI}, p. 114) and it provides a more explicit relation .

\begin{lemma} 
\label{lemmabradcond}
Let $0<a<b$ and $\epsilon>0$.
If for some ${\cal D}\in {\cal B}(\T_{W'})$ it holds
\begin{eqnarray*}
\forall {\cal E}\in {\cal B}(\T_{W^c}),\ \forall \Aa\in {\cal B}(\T_{W'}) : 
\Aa\subseteq {\cal D},\ \bP ( \Aa)>0 
\Rightarrow |\bP( {\cal E}|\Aa)-\bP({\cal E}) |< \epsilon 
\end{eqnarray*} 
then $\beta (a,b) < \epsilon \, \bP ({\cal D}) + \bP ({\cal D}^c)$. \hfill $\Box$
\end{lemma}

\medskip

\section{STIT tessellations and their $\beta$-mixing rates}\label{sectionSTIT}

For the first time, STIT tessellations were defined in
\cite{nw}. There was described a construction of STIT in bounded windows 
in all detail.  An alternative but equivalent construction 
was given in \cite{martnag}. 

\medskip

Here, let us recall  roughly the construction of the STIT tessellation process 
in a bounded window. First, we have to introduce some notation.

\medskip

Let $\hH$ denote the set of all hyperplanes in $\RR^\ell$.
We will use a parameterization of hyperplanes. 
Let $\|\cdot \|$ denote the Euclidean norm, $\langle \cdot , \cdot \rangle$ the inner product and
${\SaS}^{\ell -1}=\{ x\in \RR^\ell: \| x\| =1\} $ 
be the unit hypersphere in $\RR^\ell$. Define
$$
H(\alpha ,u )= \{ x\in \RR^\ell: \langle x,u \rangle =\alpha \}, 
\quad  \alpha \in [0,\infty ), u\in {\SaS}^{\ell -1},
$$
which is the hyperplane with normal direction $u$ 
and distance (in direction $u$) $\alpha$ from the origin. 
Notice that $H(0, u)=H(0, -u)$. Thus we can write 
$$
\hH = \left\{ H(\alpha ,u ):\, (\alpha ,u)\in [0,\infty ) \times {\SaS}^{\ell -1} \right\}
$$
and on $\hH$ we use the $\sigma$-algebra that is induced 
from the Borel $\sigma$-algebra on the parameter space. 
Any hyperplane generates two closed half-spaces
$$
H^-(\alpha ,u )= \{ x\in \RR^\ell: \langle x,u \rangle \leq \alpha \} 
\; \hbox{ and } \;
H^+(\alpha ,u )= \{ x\in \RR^\ell: \langle x,u \rangle \geq\alpha \} .
$$
For an hyperplane $H$ the above notions are only written by $H^-$ and $H^+$.

\medskip

Now, let $\Lambda$ be a (non-zero) measure
on the space of hyperplanes $\hH$ in $\RR^\ell$. It is assumed that
\begin{itemize}
\item $\Lambda$ is translation invariant;

\item $\Lambda$ possesses the following locally finiteness property:
\begin{equation}
\label{locfin1}
\Lambda([B])\!<\!\infty , \, \forall \hbox{ bounded sets } 
B\!\subset \!\RR^\ell, 
\hbox{ where }  [B]\!=\!\{H\!\in \!{\cal H}: H\cap B\neq \emptyset \};
\end{equation}

\item the support of $\Lambda$ is such that there is
no line in $\RR^\ell$ with the property that all the hyperplanes
of the support are parallel to it.

\end{itemize}

\medskip

The image of a non-zero, locally finite and
translation invariant measure $\Lambda$ with respect to 
this parameterization can be written as the product measure
\begin{equation}
\label{prodmeas}
\gamma \cdot \lambda \otimes \theta, 
\end{equation}
where $\gamma >0$ is a constant, $\lambda$ is (the restriction of) 
the Lebesgue measure on 
$[0,\infty)$ and $\theta$ is an even  probability measure on 
${\SaS}^{\ell -1}$ (cf, e.g. \cite{sw}, Theorem 4.4.1 and Theorem 13.2.12). 
Here $\theta$ is even means 
$\theta (A)= \theta (-A)$ for all Borel sets $A\subseteq {\SaS}^{\ell -1}$.

\medskip

The property that there is no line in $\RR^\ell$ 
such that all hyperplanes of the support of $\Lambda$ 
are parallel to it, is equivalent to the property that $\theta$ is 
not concentrated on a great subsphere of 
${\SaS}^{\ell -1}$, i.e there is no one-dimensional 
subspace $L_1$ of $\RR^\ell$ (with the orthogonal complement 
$L_1^\bot$) such that the support of $\theta$ equals 
${\mathbb G}= L_1^\bot \cap {\SaS}^{\ell -1}$. This property
allows to obtain a.s. bounded cells in the constructed tessellation, cf.
\cite{sw}, Theorem 10.3.2, which can also be applied to STIT tessellations.

\medskip

The assumptions made on $\Lambda$ imply 
$0<\Lambda([W])<\infty \, \hbox{ for every window } \, W$. 
Denote by $\Lambda_{[W]}$ the restriction of $\Lambda$ to $[W]$ and by 
$\Lambda^W=\Lambda ([W])^{-1}\Lambda_{[W]}$ 
the normalized probability measure. 

\medskip

Let us take a family of independent random variables
$\tau=(\tau_{n}: n\in \NN)$, where each $\tau_{n}$ is exponentially distributed 
with parameter 1.

\medskip

We will denote the random process of STIT tessellations in $W$ by 
$Y\wedge W=(Y_t\wedge W: t\geq 0)$.

\medskip

\begin{enumerate}
\item[(1)] Even if for $t=0$ the STIT tessellation $Y_0$ is not 
defined in $\RR^\ell$, we define $Y_0\wedge W= \{W\}$ the 
trivial tessellation for the window $W$. Its unique cell is denoted by 
$C^1=W$.

\item[(2)] Any extant cell has a random lifetime, 
and at the end of its lifetime it is divided by a random hyperplane. 
The lifetime of $W=C^1$ is $\Lambda ([W])^{-1} \tau_{1}$, 
and at that time it is divided by a random hyperplane $H_{1}$ 
with law $\Lambda^W$ into two 
cells denoted by $C^2=C^1\cap H_1^+$ and $C^3=C^1\cap H_1^-$.

\item[(3)] Now, any cell $C^i$ which is 
generated in the course of the construction has the lifetime 
$\Lambda ([C^i])^{-1}\, \tau_i$, i.e. it has an exponentially distributed 
lifetime with parameter $\Lambda ([C^i])$. At the end of its lifetime it 
is divided by a random hyperplane $H_i$ with law 
$\Lambda^{C^i}$. This random hyperplane is 
conditionally independent of all the lifetimes and all the dividing 
hyperplanes which appear before the present time.

\item[(4)] This procedure is 
performed for any extant cell independently. 
\end{enumerate}

\medskip

With this notation at each time $t>0$ the tessellation $Y_t\land W$ is
constituted by the cells $C^i$ which 'live' at time $t$.
It is easy to see that at any time a.s. at most one cell 
dies and so a.s. at most only two cells are born. 

\medskip

On every window $W$ there exists $Y\land W =(Y_t \land W: t > 0)$,
which we call a STIT tessellation process. It turns out to
be a pure jump Markov process and hence has the strong Markov property
(see \cite{brei}, Proposition 15.25). Each
marginal $Y_t\land W$ takes values in $\T\wedge W$.
Furthermore, for any $t>0$ the law of $Y_t$ is consistent with
respect to the windows, that is if $W'$ and $W$ are windows such that
$W'\subseteq W$, then $(Y_t \wedge W) \wedge W' \sim Y_t \wedge W'$, where
$\sim$ denotes the identity of distributions (for a proof see \cite{nw}).
This yields the existence of a STIT tessellation $Y_t$ 
of $\RR^\ell$ such that for all windows $W$ the law
of $Y_t \wedge W$ coincides with the law of the construction in the window.
A global construction for a STIT process was provided in \cite{mnw}.
A STIT tessellation process $Y=(Y_t: t>0)$ is a Markov process
and each marginal $Y_t$ takes values in $\T$.

\medskip

In the following for each tessellation $T\in \T\wedge W$ we denote
by $\C(T)$ the set of cells of $T$. We also put
$\xi(T)$ the number of cells of $T$, so by numbering the cells 
in $\C(T)$ we can write $\C(T)=\{ C^i(T):  i=1,..,\xi(T)\}$.
Note that $\xi(T)-1$ is the
number of hyperplanes that divide the window $W$ for tessellation $T$.
We define the function $\zeta(T)=\sum_{i=1}^{\xi(T)} \Lambda ([C^i(T)])$.

\medskip

In \cite{nw} it was shown that the STIT process $Y$ has no explosion, so 
at each time $t>0$ the number of cells $\xi_t:=\xi(Y_t \wedge W)$
of $Y_t\land W$, is finite a.s..
We can write by $\{C_t^i: i=1,...,\xi_t\}$ the set of cells 
of $Y_t\land W$, that is $C_t^i=C^i(Y_t\land W)$. So,
\begin{equation}
\label{defzeta}
\zeta(Y_t \wedge W) = \sum_{i=1}^{\xi_t} \Lambda ([C_t^i]).
\end{equation}
In what follows, this will be an important quantity, because (due to the 
memoryless property of the exponential distribution) the waiting time 
until the next division in $W$ after time $t$ is 
exponentially distributed with parameter $\zeta(Y_t \wedge W)$.

\medskip

Let $t_0>t$ be the time of 
the first division after $t$. If the cell $C_t^{i_0}$ is divided by 
the hyperplane $H_{i_0}$ into the cells $C_t^{i_0}\cap H_{i_0}^+$ and 
$C_t^{i_0}\cap H_{i_0}^-$, 
then the additivity of $\Lambda$ implies that
\begin{eqnarray}
\nonumber
\zeta (Y_{t_0} \wedge W) &=& \sum_{i=1}^{\xi_{t_0}} \Lambda ([C_t^i]) = 
\zeta (Y_t \wedge W)+ 
\left(\Lambda ([C_t^{i_0}\cap H_{i_0}^+])+\Lambda ([C_t^{i_0}\cap H_{i_0}^-])
-\Lambda ([C_t^{i_0}]\right)\\
\label{zetat}
&\ge & \sum_{i=1}^{\xi_{t}} \Lambda ([C_t^i]) = \zeta (Y_{t} \wedge W),
\end{eqnarray}
i.e. $\zeta$ is monotone in $t$.

\medskip

\noindent Examples: 
\begin{enumerate}
\item[(a)] {\em Isotropic model.} Assume that $\Lambda$ is also invariant 
w.r.t. rotations of the hyperplanes. Then the directional distribution 
$\theta$  on ${\SaS}^{\ell -1}$ is the uniform distribution, and, 
up to a constant factor, $\Lambda ([C_t^i])$ is the mean width of the 
cell $C_t^i$. In particular, if $\ell =2$ and $\gamma =2\pi$, 
then $\Lambda ([C_t^i])$ is the perimeter of $C_t^i$, and \\
$\zeta (Y_t \wedge W)= $ length of $\partial W$ + $2\times $
total length of edges of $Y_t$ in the interior of $W$. 

\item[(b)] If $\gamma =2\ell$ and the directional distribution is 
$\theta =\frac{1}{2\ell}\sum_{i=1}^{2\ell} \delta_{u_i}$ which is the 
uniform discrete distribution concentrated on the directions  
$u_1 ,\ldots u_{2\ell}$ of the $2\ell$ orthogonal half-axes, then 
all cells of $Y_t$ are $\ell$-dimensional cuboids, and 
$2^{\ell -1}\Lambda ([C_t^i])$ is the sum of the edge lengths of 
$C_t^i$. If $\ell =2$ then all cells are rectangles and
\\
$2 \zeta (Y_t \wedge W)=$ sum of the perimeters of all cells. 
\end{enumerate}

\medskip

\subsection{ Independent increments relation}
\label{indincr}

The name STIT is an abbreviation for "stochastic stability under
the operation of
iteration of tessellations". Closely related to that stability is a certain
independence of increments of the STIT process in time, a property which
will be used further for stating our results on $\beta-$mixing.

\medskip

In order to explain the operation of iteration, we number 
the cells of a tessellation $T\in \T$ in the following way.
Assign to each cell of $T$ a reference point in its interior 
(e.g. the Steiner point, see \cite{sw}, p. 613, or
another point that is a.s. uniquely defined). Order the set
of the reference points of all cells of $T$ by their
distances from the origin. For random homogeneous tessellations
this order is a.s. unique. Then number the cells of
$T$ according to this order, starting with
number 1 for the cell which contains the origin.
Thus we write $C^1(T),C^2(T),\ \ldots$ for the cells of $T$.

\medskip

For $T\in \T$ and
${\vec R}=(R^m: m\in \NN)\in \T^\NN$,
we define the tessellation $T\boxplus {\vec R}$, referred to
as the iteration of $T$ and ${\vec R}$, by
its set of cells
\begin{eqnarray*}
T\boxplus {\vec R}
&=& \{ C^k(T)\!\cap \!C^l({R}^k): \,
k\!=\!1,...;\, l\!=\!1,...;\,
\Int(C^k(T)\!\cap \!C^l({R}^k))\!\neq \!\emptyset \} \\
&=& \bigcup_{k \in \NN} {R}^k \wedge C^k(T).
\end{eqnarray*}
So, we restrict ${R}^k$ to the cell $C^k(T)$, and this is done for all
$k=1,\ldots $.
The same definition holds when the tessellation and the sequence of
tessellations are restricted to some window.

\medskip

To state the independence relation of the increments of the Markov
process $Y$ of STIT tessellations,
we fix a copy of the random process $Y$ and let
${\vec Y}'=({Y'}^m: m\in \NN)$ be a sequence of
independent copies of $Y$, all of them
being also independent of $Y$. In particular ${Y'}^m\sim Y$. For a fixed
time $s>0$, we set ${\vec Y}'_s=({Y_s'}^m: m\in \NN)$. Then, from the
construction and from the consistency property  of
$Y$ it is straightforward to see that the following property holds
\begin{equation}
\label{iterate}
Y_{t+s} \sim Y_t\boxplus {\vec Y}'_s \ \mbox{ for all }t,s>0\,.
\end{equation}
This relation was firstly stated in Lemma $2$ in \cite{nw}.
It implies $Y_{2t}\sim Y_t\boxplus {\vec Y}'_t$. 
The STIT property means that
\begin{equation}
\label{stit}
Y_{t} \sim 2(Y_t\boxplus {\vec Y}'_t) \ \mbox{ for all }t>0\,,
\end{equation}
so $Y_{t}\sim 2Y_{2t}$.
Here the multiplication with $2$ stands for the transformation
$x\mapsto 2x$, $x\in \RR^\ell$. 

\medskip

\section{Separating Hyperplanes}\label{sephyper}
In order to study spatial mixing properties of STIT tessellations 
we developed in \cite{martnag14} the concept of encapsulation of the 
window $W'$ inside a larger window $W$. For this, the $\Lambda$-measure 
values of sets of hyperplanes that separate facets of $W'$ from facets 
of $W$ play a key role.

\medskip

For two Borel sets $A,B\subset \RR ^\ell$ we denote 
$$
[A]= \{H\!\in \!{\hH}: H\cap A\not= \emptyset \}
$$
the set of all hyperplanes hitting $A$, and
$$
[A|B]=\{H\!\in \!{\hH}: \left(A\!\subseteq \! H^+ 
\!\land \!B \!\subseteq \! H^- \right) 
\lor \left( A\!\subseteq \! H^- \!\land \!B\!\subseteq \!H^+ \right)\} ,
$$ 
the set of all hyperplanes that separate $A$ and $B$. This set is 
a Borel set in $\hH$.

\medskip

Consider the windows $W'=[-a,a]^\ell$, $W=[-b,b]^\ell$ with $0<a<b$ and denote 
their $(\ell -1)$-dimensional facets by $f_i'$ and $f_i$ respectively, 
$i=1,\ldots ,2\ell$. We define them for $i=1 ,\ldots ,\ell$ as 
$$
f_i'= [-a,a]\times \ldots \times [-a,a]\times \{ a\} 
\times [-a,a]\times \ldots \times [-a,a]
$$ 
with the singleton $\{ a\}$ standing on the $i$-th position,
and 
$f_{i+\ell}'= -f_i'$ for $i=1 ,\ldots ,\ell$,
and the $f_i$ are defined as the $f_i'$ respectively, by replacing $a$ by $b$.
We will use the sets of separating hyperplanes 
\begin{equation}\label{defG}
G_i(a,b)=[f_i' | f_i], 
\quad i=1,\ldots, 2\ell  
\end{equation} 

\begin{lemma}
\label{disjointG}
If $0<a<b$ then the sets of hyperplanes 
$G_i(a,b)=[f_i' | f_i], \, i=1,\ldots, 2\ell,$ are all nonempty, 
and they are pairwise disjoint. Furthermore for $r>1$ 
and $ i=1,\ldots , 2\ell$:
\begin{equation}
\label{suplinearity}
\Lambda (G_i( a, r\, b)) \geq  
r \Lambda (G_i(a,b)), 
\end{equation} 
and
\begin{equation}
\label{linearity}
\Lambda (G_i( r\, a, r\, b)) =  r\, \Lambda (G_i(a,b)).
\end{equation}
\hfill $\Box$
\end{lemma}

For the rest of this paper we need the following
additional assumption on $\Lambda$: For all $0<a<b$ and the 
windows $W'=[-a,a]^\ell$, $W=[-b,b]^\ell$ we assume that 
\begin{equation}
\label{assLambda}
\Lambda (G_i(a,b))>0 \quad \mbox{ for all } i=1,\ldots , 2\ell .
\end{equation}

\medskip

This condition is fulfilled, e.g., if $\Lambda$ is rotation invariant 
(isotropic STIT model) or if it is concentrated on the sets of 
hyperplanes that are orthogonal to coordinate axes (see Remark 1 (iii) 
in \cite{martnag14}). The condition is not fulfilled, e.g., when 
in $\RR^2$ the measure $\Lambda$ is concentrated on the two directions 
that are parallel to the diagonals of $W'$ and $W$.

\medskip

Note that due to the translation invariance of $\Lambda$ and the symmetry 
of $W'$ and $W$ we have $\Lambda (G_i(a,b)) = \Lambda (G_{i+\ell}(a,b))$ 
for all $i=1,\ldots ,\ell $.
Let us define
\begin{equation}
\label{defL}
L(a,b)=\min \{ \Lambda (G_i(a,b)),\, i=1,\ldots , 2\ell \} .
\end{equation}
Obviously, $L(a,b)$ has the properties analogous to (\ref{suplinearity}) and 
(\ref{linearity}) as well.

\section{Main Results}

Now we are prepared to formulate our main results. But beforehand 
let us state some notation precisely. The $\beta(a,b)$ defined in (\ref{defbeta}) 
obviously depends on the probability measure $\bP$. In our context where 
we study the mixing for the STIT tessellation process at some fixed 
time $t>0$, the $\beta(a,b)$ coefficient depends on 
$\PP(Y_{t}\in dT)$, which is the distribution of $Y_t$. Hence it depends 
on $t>0$, but as time is fixed we shall not indicate explicitly the dependence of 
$\beta(a,b)$ on $t$.

\begin{theorem}
\label{TheorSTITbeta}
Let $Y$ be the STIT tessellation process determined by the hyperplane measure 
$\Lambda$ satisfying (\ref{assLambda}). Further, let $\zeta (Y_t \wedge W')$ 
be as defined in (\ref{defzeta}) and $L$ as defined in (\ref{defL}) for $0<a<b$, 
$W'=[-a,a]^\ell \subset W= [-b,b]^\ell$. Then for a fixed $t>0$ and all 
$0<s<t,\ M>0$  we have
\begin{eqnarray*}
\beta (a,b)& < & \PP ( \zeta (Y_t\wedge W') \geq M)\\
&{}& \, + \, \PP ( \zeta (Y_t\wedge W') < M)\cdot
\left[ 1- {\rm e}^{-s \Lambda ([W'])} 
\left(1- {\rm e}^{-s L(a,b)} \right)^{2 \ell} {\rm e}^{-s M} \right.\\
&{}& \;\;\;\; + \max \left\{ \left( {\rm e}^{s M}- 1 \right) ;  
2-{\rm e}^{-s M}  - {\rm e}^{-s \Lambda ([W'])}   
\left(1\!- \!{\rm e}^{-s L(a,b)} \right)^{2 \ell} \right\} \! \Biggr].
\end{eqnarray*} \hfill $\Box$
\end{theorem}

If the purpose is to give a bound for the order of decay of $\beta$, 
instead of the $\max$ one can consider the sum of the two nonnegative items. 
Hence one can also use the bound
\begin{eqnarray*}\label{simplupper}
\beta (a,b) &<& \PP ( \zeta (Y_t\wedge W')\! \geq \!M)  + \, 
\PP (\zeta (Y_t\wedge W')\!< \!M)\times \\
&{}& \quad\quad 
\times \left[\! 2 + {\rm e}^{s M}- {\rm e}^{-s M} -(1+\!{\rm e}^{-s M}) 
{\rm e}^{-s \Lambda ([W'])}\left(1- {\rm e}^{-s L(a,b)}\right)^{2 \ell} 
\right].
\end{eqnarray*}
This upper bound can now be minimized by choosing 
appropriate $M>0$ and $0<s<t$. 
A (rough) upper bound for $\PP ( \zeta (Y_t\wedge W') \geq M)$ 
can be derived with the Chebyshev inequality (see the proof of Theorem \ref{decaypoly}).

\bigskip

The following result states that for STIT tessellations the condition 
(\ref{bradcond}) is satisfied. 

\medskip

\begin{theorem}
\label{eqivbetamix}
Let $Y$ be a STIT 
tessellation process determined by the hyperplane measure 
$\Lambda$ with (\ref{assLambda}). Then  for $0<a<b$, 
$W'=[-a,a]^\ell \subset W= [-b,b]^\ell$ 
\begin{eqnarray*}
&{}&\forall a>0,\ \forall t>0, \ \forall \epsilon>0,\ 
\exists {\cal D}\in {\cal B}(\T_{W'}) ,\ 
\PP (Y_t \in {\cal D})>1-\epsilon ,\ \exists b>a,\\ 
&{}&\forall {\cal E}\in {\cal B}(\T_{W^c}),\ 
\forall \Aa\in {\cal B}(\T_{W'}) : \\
&{}& \Aa\subseteq {\cal D},\ \PP (Y_t \in \Aa)>0
\Rightarrow |\PP(Y_t\in {\cal E}\, | \, Y_t\in \Aa)-
\PP(Y_t\in {\cal E}) |
< \epsilon .
\end{eqnarray*} \hfill $\Box$
\end{theorem}

Together with Theorem \ref{bradequiv} it immediately implies the first part of the
following result.

\begin{theorem}
\label{decaypoly}
Under the assumptions of Theorem \ref{TheorSTITbeta} we have
$\lim_{b\to \infty} \beta (a,b) =0$ for all $a>0$. 

\medskip

\noindent Moreover for all $\eta\in (0,1)$ 
there exists a constant $\kappa=\kappa(t,a,\eta)<\infty$
such that $\beta (a,b)\le \kappa \, b^{-\eta}$ for all $b>a$.
\hfill $\Box$
\end{theorem}

\section{Proofs}\label{sectionproofs}

\noindent {\bf Proof of Lemma \ref{disjointG}:} 

\noindent
(i) Assume that $u_i \in {\SaS}^{\ell -1}$ is the normal vector 
for $f_i', f_i$. Then $\emptyset \not= \{ H(\alpha ,u_i):\, 
a\leq \alpha \leq b \} \subset G_i(a,b) $.

\medskip

\noindent (ii) In order to show that the $G_i(a,b)$, $i=1,\ldots , 2\ell$, 
are pairwise disjoint, we consider the cones 
$$
K_i= \{ \lambda x:\, \lambda \geq 0,\ x\in f_i\}
$$
and we show that $H(\alpha ,u)\in G_i(a,b)$ implies 
that $u \in int (K_i)$, the interior of $K_i$.
For simplicity, put $i=1$ and write
$$
f_1' = a\cdot conv( \{ (1,(-1)^{\beta_2},\ldots ,(-1)^{\beta_\ell}): 
\beta_j\in \{ 0,1\} \} )
$$
where $conv$ denotes the convex hull of a set.
If $H(\alpha ,u)\in G_1(a,b)$, then in particular it 
separates all pairs of vertices with one vertex in $f_1'$ 
and the other vertex in $f_1$. Hence $H(\alpha ,u)$ intersects 
all the 1-dimensional edges of the cone $K_1$, i.e. for all 
$j=2,\ldots ,\ell$ and the corresponding (to the $j$-th edge of $K_i$) 
choice of $\beta_2,\ldots ,\beta_\ell$ there exists a $\lambda_j >0$ such that
$$
\lambda_j \cdot (1, (-1)^{\beta_2},\ldots ,(-1)^{\beta_\ell}) \in H(\alpha ,u) .
$$
Thus the definition of the hyperplanes implies for $u=(u^1,\ldots , u^\ell)$ that
$$
 u^1 + \sum_{i=2}^\ell (-1)^{\beta_i}\, u^i =\frac{\alpha}{\lambda_j } 
 > 0.
$$
Because this must be satisfied for all edges of $K_1$, i.e. for all choices of
$(\beta_2,\ldots ,\beta_\ell)\in \{ 0,1\}^{\ell -1}$, the vector $u$ must 
fulfill the condition
$$
u^1 - \sum_{i=2}^\ell |u^i| > 0.
$$
Furthermore, $\| u\| =1$ implies $u^1=\sqrt{1-\sum_{i=2}^\ell  (u^i)^2}$ 
and hence
$$
1 > 2 \left(  \sum_{i=2}^\ell  (u^i)^2  +  \sum_{2\leq i<k\leq \ell}  
|u^i \cdot u^k| \right) .
$$
This yields $\sum_{i=2}^\ell  (u^i)^2<\frac{1}{2}$ and hence $u^1 >\frac{1}{2}\sqrt{2}$ and $u^i<\frac{1}{2}\sqrt{2}$ 
for all $i=2, \ldots ,\ell$ which implies that $u\in int (K_1)$.

\medskip

\noindent (iii) For a nonempty closed convex set $A\subset \RR ^\ell$ the support 
function $h(A,\cdot )$ is defined  by
$$
h(A,u)=\sup \{ \langle x,u \rangle :\, x\in A \}, \quad  u\in {\SaS}^{\ell -1} .
$$
Notice that
$$
\inf \{ \langle x,u\rangle  :\, x\in A \}=
-\sup \{ \langle x,-u\rangle  :\, x\in A \} = - h(A,-u),
$$
and for homothets
\begin{equation} 
\label{homothet}
h(r\cdot A, u) =\sup \{ \langle x,u\rangle  :\, x\in r\cdot A \} = r \, 
h(A,u) , \quad  u\in {\SaS}^{\ell -1}, r>0 .
\end{equation}
 To show (\ref{suplinearity}) and (\ref{linearity}), notice 
that for a given pair $f_i',\, f_i$ and direction $u \in int(K_i)$ a 
hyperplane $H(\alpha , u)\in [f_i'|f_i]$, if and only if 
$\alpha \geq h(f_i',u)$ and 
$\alpha \leq \inf \{ \langle x,u \rangle :\, x\in f_i \}=
-\sup \{ \langle x,-u\rangle  :\, x\in f_i \}=-h(f_i,-u)$. Hence, (\ref{prodmeas}) yields
$$
\Lambda ([f_i'|f_i]) =\gamma \int_{K_i} \left[ -h(f_i,-u)- h(f_i',u) \right]_+ 
\theta (du) ,
$$
where $[\cdot ]_+$ denotes the nonnegative part.
Using (\ref{homothet}) immediately completes the proof of (\ref{linearity}).
Furthermore, for $u\in K_i$, because $h(f_i',u)>0$ , and $r>1$,
\begin{eqnarray*}
&{}& \left[ -h(rf_i,-u)- h(f_i',u) \right]_+ \\ 
& = & \left[ -r h(f_i,-u) -r h(f_i',u) + r h(f_i',u)- h(f_i',u) \right]_+ \\ 
& \geq & \left[ -r h(f_i,-u) -r h(f_i',u)\right]_+
\end{eqnarray*}
and this yields (\ref{suplinearity}). \hfill $\Box$
\bigskip

Now we prepare the proof of Theorem \ref{TheorSTITbeta} by several 
lemmas. 

\medskip

In the sequel $Y$ is the STIT process and $W'$ is a window. Let $t>0$ and 
$s\in [0,t)$. 

\medskip

Let $\Aa\in {\cal B}(\T_{W'})$ be a measurable set such that
$\PP({Y_t}\in \Aa)>0$.
It is easy to see that this implies $\PP({Y_{t'}}\in \Aa)>0$
for all $t'>0$. 

\medskip

Let us first provide a result only concerning the process on a window $W'$.

\medskip

We define
\begin{equation}
\label{not1}
\chi^\Aa(t,s;t')=\PP( Y_t\wedge W'=Y_{t-s}\wedge W' \, | \, Y_{t'}\in \Aa)
\;\hbox{ for } t'\in \{t-s,t\}.
\end{equation}
Note that the following equality of events is satisfied
\begin{equation}
\label{not2}
\{Y\wedge W' \hbox{ has no jumps in } (t-s,t]\}=
\{Y_t\wedge W'=Y_{t-s}\wedge W'\}.
\end{equation}

\begin{lemma}
\label{lemmajump}
We have
\begin{equation*}
\chi^\Aa(t,s;t) \leq   \frac{\PP(Y_{t-s}\in \Aa)}{\PP(Y_t\in \Aa)}
\leq \chi^\Aa(t,s;t\!-\!s)^{-1}.
\end{equation*}
\end{lemma}

\noindent
{\bf Proof:}
\begin{eqnarray*}
\PP(Y_{t-s}\in \Aa)
&= & \PP(Y_{t-s}\in \Aa,\ Y\wedge W' \hbox{ has no jump in } (t-s,t])  \\
&{}& \;\, + \, \PP(Y_{t-s}\in \Aa,\ Y\wedge W' \hbox{ has a jump in } (t-s,t]) \\
&\geq & \PP(Y_t\in \Aa, \ Y\wedge W' \hbox{ has no jump in } (t-s,t]).
\end{eqnarray*}
Dividing both sides by $\PP(Y_t\in \Aa)$ provides the first inequality.
On the other hand,
$$
\PP(Y_t\in \Aa) \geq \PP(Y_{t-s}\in \Aa, \ Y\wedge W'
\hbox{ has no jump in } (t-s,t])
$$
which yields
$\PP(Y_t\in \Aa)/\PP(Y_{t-s}\in \Aa) \geq \chi^\Aa(t,s;t\!-\!s)$,
and taking the reciprocal on both sides completes the proof. 
\hfill $\Box$

\begin{remark}
This lemma can be stated for the process $Y\wedge W'$, 
or more generally, for any  pure jump process $Z=(Z_t: t\ge 0)$ 
taking values in some set $V$, and $\Aa\subseteq V$ a measurable set 
such that $\PP(Z_{t'}\in \Aa )>0$ for $t'\in \{t-s,t\}$. Analogously as 
in (\ref{not1}) we put
$\chi^\Aa(t,s;t')=
\PP( Z \hbox{ has no jump in } (t-s,t] \, | \, Z_{t'}\in \Aa)$ for 
$t'\in \{t-s,t\}$.
But in this case the condition (\ref{not2}) does not necessarily hold.
Under these conditions the same proof as the one of Lemma \ref{lemmajump}
gives 
$\chi^\Aa(t,s;t) \leq   \PP(Z_{t-s}\in \Aa)/\PP(Z_t\in \Aa) 
\leq \chi^\Aa(t,s;t\!-\!s)^{-1}$. \hfill $\Box$
\end{remark}

\medskip

From now on, we fix two windows $W', W$, such that $W'\subset W$.
Let us give some results concerning the encapsulation time of $W'$ inside $W$. 

\medskip

In \cite{martnag14} we introduced the concept of encapsulation of $W'$ 
inside $W$, which means that there is a state of the process $Y$ 
(or, equivalently, 
of $Y\wedge W$) such that all facets of $W'$ are separated from the 
facets of $W$ by facets of the tessellation {\em before} the interior 
of $W'$ is divided by a facet of the tessellation. Formally, denoting 
the cell of $Y_t$ that contains the origin (the 0-cell) by $C_t^1$, 
we define the encapsulation time as
$$
S(W',W) =\inf \{ t>0:\, W'\subseteq C_t^1 \subset int(W) \} ,
$$  
and we put $\inf \emptyset =\infty$.
Notice that $Y_{t'}\wedge W' = W'$ describes the event that 
until time $t'$ there is no division of $W'$ in the process $Y$ 
(or, equivalently, in the process $Y\wedge W$).

\begin{lemma}
\label{lemmacondencaps}
For all ${\cal E}\in {\cal B}(\T_{W^c})$ we have
\begin{eqnarray}
\label{eqsolita}
&{}  &\PP(Y_t\in {\cal E}, S(W',W)\!<\!s, Y_{s}\wedge W' \!=\! W' \, | \, Y_t\!\in \!\Aa) \\
&{}= & \PP(Y_t\in {\cal E}, S(W',W)\!<\!s, Y_s\wedge W' \!=\! W')
\frac{\PP( Y_{t-s}\in \Aa)}{\PP( Y_t\!\in \!\Aa)}. \nonumber
\end{eqnarray}
\end{lemma}

\noindent
{\bf Proof:}
We will use the notation and the result on the independent increments 
relation stated in Section \ref{indincr}.

\medskip

By $Y_s\boxplus (Y'^m_{t-s}: m \geq 2)$ we mean
that the tessellations $Y'^m_{t-s}$ are nested only
into the cells $C^m_s$ of $Y_s$ with $m \geq 2$, and not
into the $0$-cell $C^1_s$. Below we use
the independence of the random variables
$Y_s$, $Y'^m_{t-s}$, $m \geq 1$, the implication
$$
(S(W',W)\!<\!s, Y_{s}\wedge W' = W) \Rightarrow
(W'\subseteq C^1_s \subset \Int W)
$$
and ${\cal E}\in {\cal B}(\T_{W^c})$, to get:
\begin{eqnarray*}
&{} & \PP(Y_t\in {\cal E},\, S(W',W)\!<\!s,
Y_s\wedge W' \!= \!W' \, | \, Y_t\in \Aa) \\
&{}= & \frac{1}{ \PP( Y_t\in \Aa)} \PP(Y_t\in {\cal E},\, S(W',W)\!<\!s,\,
Y_{s}\wedge W' \!= \!W',\,  Y_t\in \Aa) \\
&{}= & \frac{1}{ \PP( Y_t\in \Aa)} \PP(Y_{s}\boxplus
{\vec Y}'_{t-s}\in {\cal E},\, S(W',W)\!<\!s,\, Y_{s}\wedge W' \!= \!W',\,
Y'^1_{t-s}\in \Aa) \\
&{}= & \frac{1}{ \PP( Y_t\in \Aa)} \PP(Y_{s}\boxplus
(Y'^m_{t-s}: m \geq 2)\in {\cal E},\,
S(W',W)\!<\!s,\, Y_{s}\wedge W' \!= \!W',\,  Y'^1_{t-s}\in \Aa) \\
&{}= & \frac{1}{ \PP( Y_t\in \Aa)}
\PP(Y_t\in {\cal E},\, S(W',W)\!<\!s,\,
Y_{s}\wedge W' \!= \!W')\cdot \PP(Y'^1_{t-s}\in \Aa)\\
&{}= & \frac{\PP(Y_{t-s}\in \Aa)}{ \PP( Y_t\in \Aa)}\,
\PP(Y_t\in {\cal E},\, S(W',W)\!<\!s,\, Y_{s}\wedge W' \!= \!W') .
\end{eqnarray*} 
\hfill $\Box$

\medskip

Equation (\ref{eqsolita}), will be used
combined with the bounds on $\PP( Y_{t-s}\in \Aa)/\PP( Y_t\in \Aa)$
provided by Lemma \ref{lemmajump}.

\medskip

\begin{lemma}
\label{lemmadiffencaps} 
For ${\cal E}\in {\cal B}(\T_{W^c})$ we have
\begin{eqnarray*}
&{}  & |\PP(Y_t\in {\cal E},\, S(W',W)\!<\!s, 
Y_{s}\wedge W' \!= \!W' \, | \, Y_t\in \Aa)  - 
\PP(Y_t\in {\cal E})| \\
&\leq & \max \Biggl\{ \left(\chi^\Aa(t,s;t\!-\!s)^{-1} - 1 \right) 
\cdot \PP(Y_t\in {\cal E}) ; \ 
(1- \chi^\Aa(t,s;t))\cdot  \PP(Y_t\in {\cal E}) \\
&{}& \qquad\qquad  + \, \chi^\Aa(t,s;t) \cdot  
(1- \PP(S(W',W)\!<\!s,\, Y_{s}\wedge W' \!= \!W'))  \Biggr\} .
\end{eqnarray*}
\end{lemma}

\noindent {\bf Proof:}
We have the inequality
$$
\PP(Y_t\in {\cal E}) -(1- \PP(S(W',W)\!<\!s,\, Y_{s}\wedge W' \!= \!W'))
\leq \PP(Y_t\in {\cal E}, S(W',W)\!<\!s,\, Y_{s}\wedge W' \!= \!W').
$$
Then, from Lemmas \ref{lemmajump} and \ref{lemmacondencaps} we get
\begin{eqnarray*}
&{}& \chi^\Aa(t,s; t)  
(\PP(Y_t\in {\cal E}) -(1- \PP(S(W',W)\!<\!s,\, Y_{s}\wedge W' \!= \!W')))\\
\leq &{}& \chi^\Aa(t,s; t)  
\PP(Y_t\in {\cal E},\, S(W',W)\!<\!s,\, Y_{s}\wedge W' \!= \!W')\\
\leq &{}& \PP(Y_t\in {\cal E},\, S(W',W)\!<\!s, 
Y_{s}\wedge W' \!= \!W' \, | \, Y_t\in \Aa) \\
\leq &{}& \PP(Y_t\in {\cal E},\, S(W',W)\!<\!s,\, Y_{s}\wedge W' \!= \!W')
\chi^\Aa(t,s; t\!-\!s)^{-1}\\
\leq &{}& \PP(Y_t\in {\cal E})\chi^\Aa(t,s; t\!-\!s)^{-1}.
\end{eqnarray*}

Subtraction of $\PP(Y_t\in {\cal E})$ in all lines of these inequalities 
leads to
\begin{eqnarray*}
&{} & (\chi^\Aa(t,s; t) -1) \PP(Y_t\in {\cal E})
- \chi^\Aa(t,s; t)  
(1- \PP(S(W',W)\!<\!s,\, Y_{s}\wedge W' \!= \!W')) \\
\leq &{}& \PP(Y_t\in {\cal E},\, S(W',W)\!<\!s, 
Y_s\wedge W' \!= \!W' \, | \, Y_t\in \Aa) - \PP(Y_t\in {\cal E}) \\
\leq &{}& \left( \chi^\Aa(t,s; t\!-\!s)^{-1} - 1 \right) 
\PP(Y_t\in {\cal E})\,,
\end{eqnarray*}
and since the term in the first line is less or equal zero, we obtain 
the result for the absolute value. 
\hfill $\Box$

\medskip

\begin{lemma}
\label{cordelta}
Let $Y$ be the STIT process and ${\cal D}\in {\cal B}(\T_{W'})$. Consider 
the set,
$$
\U=\{\Aa\in {\cal B}(\T_{W'}):  \Aa\subseteq {\cal D}, \PP ( Y_t\in \Aa)>0\},
$$
and so $\U\times {\cal B}(\T_{W^c})=\{(\Aa,{\cal E}): \Aa\in \U, 
{\cal E}\in {\cal B}(\T_{W^c})\}$. Then
\begin{eqnarray*}
&{} &\sup\limits_{\U\times {\cal B}(\T_{W^c})} 
|\PP(Y_t\in {\cal E}\, | \, Y_t\in \Aa)-\PP(Y_t\in {\cal E}) |  \\
&{}\leq & \sup\limits_{\Aa\in \U} 
(1- \PP( S(W',W)\!<\!s, Y_{s}\wedge W' \!= \!W' \, | \, Y_t\in \Aa) ) \\
&{}& \;\; + \sup\limits_{\U\times {\cal B}(\T_{W})} 
\max \Biggr\{ 
\left(\chi^\Aa(t,s; t\!-\!s)^{-1} - 1 \right)\PP(Y_t\in {\cal E});  
  \\
&{}& \;\;\;\; (1- \chi^\Aa(t,s; t) ) 
\PP(Y_t\in {\cal E}) + \chi^\Aa(t,s; t) 
(1- \PP(S(W',W)\!<\!s, Y_{s}\wedge W' \!=\! W'))  \Biggr\}.
\end{eqnarray*}
\end{lemma}

\noindent {\bf Proof:}
\begin{eqnarray*}
&{}  &\sup\limits_{\U\times {\cal B}(\T_{W})} 
|\PP(Y_t\in {\cal E}\, | \, Y_t\in \Aa)-\PP(Y_t\in {\cal E}) |  \\
&{}\leq & \sup\limits_{\U\times {\cal B}(\T_{W})} 
|\PP(Y_t\in {\cal E}\, | \, Y_t\in \Aa)-\PP(Y_t\in {\cal E},\, 
S(W',W)\!<\!s, Y_{s}\wedge W' \!= \!W' \, | \, Y_t\in \Aa) | \\
&{} & \; + \sup\limits_{\U\times {\cal B}(\T_{W})} 
|\PP(Y_t\in {\cal E},\, S(W',W)\!<\!s, 
Y_{s}\wedge W' \!= \!W' \, | \, Y_t\in \Aa)-
\PP(Y_t\in {\cal E}) |.
\end{eqnarray*}
Now we use the general inequality $\PP(A \, | \, B)-\PP(A\cap C \, | \, B)
\le 1-\PP(C \, | \, B)$ for $A=\{Y_t\in {\cal E}\}$, $B=\{Y_t\in \Aa\}$ 
and $C=\{S(W',W)\!<\!s, Y_{s}\wedge W' \!= \!W'\}$, and 
then Lemma \ref{lemmadiffencaps} to conclude
the proof. 
\hfill $\Box$

\medskip

In \cite{martnag14}, last part of the proof of Lemma 4, it was shown that for STIT processes
\begin{eqnarray}
\nonumber
\PP(S(W',W)\!<\!s,\ Y_{s}\wedge W' \!= \!W') &\geq&
{\rm e}^{-s \Lambda ([W'])}
\prod_{i=1}^{2 \ell} \left(1- {\rm e}^{-s \Lambda (G_i(a,b))} \right) \\
\label{encapstime}
&\geq&  {\rm e}^{-s \Lambda ([W'])} \,  
\left(1- {\rm e}^{-s L(a,b)} \right)^{2 \ell}.
\end{eqnarray}

In the next result we make use of 
the notion $\zeta(T\wedge W')$ defined in (\ref{defzeta}).

\begin{lemma}
\label{lemmanojumpcondA}
For any $M>0$ and ${\cal D}=  \{ T \in \T : \zeta(T\wedge W') < M \}$  
we have for all $\Aa\in {\cal B}(\T_{W'}),\, \Aa\subseteq {\cal D},\ 
\PP(Y_t\in \Aa)>0$
\begin{equation}
\label{casos}
\chi^\Aa(t,s;t\!-\!s)\geq {\rm e}^{-s M}\;\; \hbox{ and } \;\; 
\chi^\Aa(t,s;t)\geq {\rm e}^{-s M}.
\end{equation}
\end{lemma}

\noindent {\bf Proof:}
As mentioned after (\ref{defzeta}), the waiting time after $t-s$ until the 
next jump of the process $Y\wedge W'$ (where there is a division in 
$Y_{t-s}\wedge W'$), is exponentially distributed with parameter 
$\zeta(Y_{t-s}\wedge W')$. By conditioning on $\sigma(Y_u: u\le t-s)$ 
and the Markov property we get,
\begin{eqnarray*}
\chi^\Aa(t,s;t\!-\!s)\PP(Y_{t-s}\in \Aa)&=&\PP(Y_{t-s}\wedge W'=Y_t\wedge W', 
Y_{t-s}\in \Aa)\\
&=&
\int_\Aa e^{-s\zeta(T\wedge W')} \PP(Y_{t-s}\in dT).
\end{eqnarray*}
(As already said $\PP(Y_{t-s}\in dT)$ the distribution of $Y_{t-s}$ on $\T$).
By definition $\zeta(T\wedge W') < M$ when $T\in {\cal D}$,  in particular
when $T\in \Aa$. Then,
$$
\chi^\Aa(t,s;t\!-\!s)\PP(Y_{t-s}\in \Aa)\ge e^{-sM}\PP(Y_{t-s}\in \Aa).
$$
The first relation in (\ref{casos}) follows.

\medskip

Now,
\begin{eqnarray*}
\chi^\Aa(t,s;t)\PP(Y_{t}\in \Aa)&=&\PP(Y_{t-s}\wedge W'=Y_t\wedge W', Y_t\in \Aa)\\
&=& \int_{\{Y_t\in \Aa\}}{\bf 1}_{\{Y_{t-s}\wedge W'=Y_t\wedge W'\}} d\PP(Y)\\
&=&\int_{\{Y_t\in \Aa\}} \EE({\bf 1}_{\{Y_{t-s}\wedge W'=Y_t\wedge W'\}}\, | \, Y_t)
d\PP(Y)\\
&=&\int_{\{Y_t\in \Aa\}} e^{-s\zeta(Y_t\wedge W')} d\PP(Y)\\
&= &\int_\Aa e^{-s\zeta(T\wedge W')} P(Y_t\in dT)
\end{eqnarray*}
(We have denoted by $\PP(dY)$ the distribution of $Y$ on the space of 
trajectories). Again we use that for all $T\in \Aa\subseteq {\cal D}$ we have 
$\zeta(T\wedge W') < M$, to conclude that
$$
\chi^\Aa(t,s;t)\PP(Y_{t}\in \Aa)\ge e^{-sM}\PP(Y_t\in \Aa).
$$
The result follows.
\hfill $\Box$

\medskip

\noindent {\bf Proof of Theorem \ref{TheorSTITbeta}:}

\noindent Applying Lemma \ref{cordelta} and using Lemma \ref{lemmacondencaps}, 
(\ref{encapstime}) and Lemma \ref{lemmanojumpcondA} yields
\begin{eqnarray}
\label{boundsup}
&{} &\sup\limits_{({\cal E}, \Aa)\in \U} 
|\PP(Y_t\in {\cal E}\, | \, Y_t\in \Aa)-\PP(Y_t\in {\cal E}) |  \nonumber \\
&{}\leq &   1-  {\rm e}^{-s \Lambda ([W'])} \,  
\left(1- {\rm e}^{-s L(a,b)} \right)^{2 \ell} 
\, {\rm e}^{-s M}  \\
&{} & + \max \left\{  {\rm e}^{s M} - 1  ;\  1- {\rm e}^{-s M}  \,   + \,  
1- {\rm e}^{-s \Lambda ([W'])} \,  
\left(1- {\rm e}^{-s L(a,b)} \right)^{2 \ell}  \right\}. 
\nonumber 
\end{eqnarray}
\hfill $\Box$

\medskip

Before proving Theorem \ref{eqivbetamix} we show an intermediate result,
there we use the notation $\zeta (Y_t\wedge W')$ introduced in (\ref{defzeta}).

\begin{lemma}
\label{setD}
For all $a>0$, $t>0$ and $\epsilon>0$ there exists $M(a,t,\epsilon )>0$
such that $\PP (\{  \zeta (Y_t\wedge W') < M(a,t,\epsilon)\} ) > 1-\epsilon$.
\end{lemma} 

\noindent {\bf Proof:}
Because the process $(Y_t\wedge W': t>0)$ has no explosion 
(this was shown in \cite{nw}), equation (\ref{zetat}) implies that  
$\zeta (Y_t\wedge W')$ is a.s. finite for all $t>0$.
Therefore, 
$$
\PP \left( \bigcup_{n=1}^{\infty} \{  \zeta(Y_t\wedge W')< n \} \right) = 
\lim_{n\to \infty }\PP (\zeta (Y_t\wedge W')< n \} )= 1,
$$
which proves the assertion. 
\hfill $\Box$

\medskip

\noindent {\bf Proof of Theorem \ref{eqivbetamix}:}

\noindent We fix $a,\, t,\, \epsilon>0$, $W'=[-a,a]^\ell$ and 
$M=M(a,t,\epsilon)$. Then choose $s\in (0,t)$ such that 
the following three conditions are satisfied:
\begin{eqnarray*}
&{}& {\rm e}^{s M} - 1 < \frac{\epsilon}{2}\,;\\
&{}&{\rm e}^{-s M} > \left(1-\frac{\epsilon}{2}\right)^{1/3}\,;\\
&{}& {\rm e}^{-s \Lambda ([W'])}> \left(1-\frac{\epsilon}{2}\right)^{1/3}\,.
\end{eqnarray*}
The existence of such an $s$ is obvious.
Now, for such an $s$, from  $\Lambda$ satisfying (\ref{assLambda}) and 
(\ref{suplinearity}) it can be chosen $b>a$ such that $L(a,b)$ is as big to 
have $\left( 1-{\rm e}^{-s L(a,b)} \right)^{2 \ell}> 
\left( 1-\frac{\epsilon}{2}\right)^{1/3}$.
Now, plugging this into (\ref{boundsup}), we obtain Theorem \ref{eqivbetamix}. 
\hfill $\Box$

\bigskip

\noindent {\bf Proof of Theorem \ref{decaypoly}:}
We must only show the second assertion on the power of decay.

\medskip

As stated, the jump rate of $\xi(Y_t \wedge W')$ is 
bounded by $q\xi(Y_t \wedge W')$ where $q=\Lambda([W'])$. In fact,
the jump rate of $\xi(Y_t \wedge W')$ is
$\zeta(Y_t \wedge W')$, and (\ref{defzeta}) implies that this 
is bounded by 
$$
\zeta(Y_t \wedge W') = \sum_{i=1}^{\xi(Y_t \wedge W')} \Lambda ([C_t^i])
\le q \xi(Y_t \wedge W').
$$
Then, the process $\xi(Y\wedge W')=(\xi(Y_t\wedge W'): t\ge 0)$ 
is stochastically dominated  by a continuous time birth chain 
$B=(B_t: t\ge 0)$ starting from $B_0=1$, 
with linear rates $(qn: n \ge 1)$. That is, 
$B$ satisfies, $\PP(B_{t+h}=n+1 \, | \, B(t)=n)=qn+o(h)$. 
(As usual $o(h)$ means $o(h)/h\to 0$ as $h\to 0^+$).
Then the process $q^{-1}\zeta(Y\wedge W')=(q^{-1}\zeta(Y_t\wedge W'): t\ge 0)$
is dominated by $B$.
This stochastic domination means
$\PP(q^{-1}\zeta(Y_t\wedge W')\ge M)\le \PP(B_t\ge M)$. 

\medskip

But the linear birth process has finite moments of all orders, in fact 
$\EE(B_t^r)=e^{-qt}\sum_{l\ge 1}  l^r (1-e^{-qt})^{l-1} <\infty$ for all 
$r\ge 1$ (see \cite{grimmett} Exercise $6.8.19$, p. $59$).
Hence, the stochastic domination of $q^{-1}\zeta(Y\wedge W')$ by $B$, 
implies that $\EE(\zeta(Y_t\wedge W')^r)<\infty$ for all $r\ge 1$.
(We note that at least in the isotropic case, 
$\EE(\zeta(Y_t\wedge W'))$ is known, and a formula for 
Var$(\zeta(Y_t\wedge W'))$ is given in \cite{schthaeLIM3}).

\medskip
  
Applying the Markov inequality, for all $r>0$ and all $M>0$ we obtain 
$$
\PP (\zeta(Y_t\wedge W')  \ge  M)\le   
\frac{\EE(\zeta(Y_t\wedge W')^r)}{M^r} ,
$$
Now put $M=b^v$ with $0<v$. Hence for every $r>0$ we have that for
$\kappa_1(t,a,r)=\EE(\zeta(Y_t\wedge W')^r)$,
$$
\PP (\zeta(Y_t\wedge W')  \ge  M)\le \kappa_1(t,a,r)b^{-vr}. 
$$
(The dependence on $a$ is because $W'=[-a,a]^\ell$).
Now put $s=b^{-u}$. with $u>0$. Then
$$
{\rm e}^{-s\Lambda([W'])}=1-\Lambda([W'])b^{-u}+o(b^{-u}).
$$
Now we assume the condition $0<v<u$ which guarantees
that $sM\to 0$ as $b\to \infty$ and also it implies
$$
{\rm e}^{s M}- {\rm e}^{-s M}=\frac{1}{2} (sM)^2+o((sM)^2)=
\frac{1}{2} b^{-2(u-v)}+o(b^{-2(u-v)}).
$$ 
Fix $b_0>0$ and let $L_0(a)=\frac{L(a,b_0)}{b_0}$ which is strictly positive.
From Lemma \ref{disjointG} we have for all $b>b_0$, 
$$
L(a,b)=L(a,(b/b_0)b_0)\ge b L_0(a).
$$
Now also take the constraint $u<1$ to have
$sL(a,b)\ge b^{1-u} L_0(a)$. Therefore 
$$
{\rm e}^{-sL(a,b)}\le {\rm e}^{-b^{1-u} L_0(a)}=({\rm e}^{b^{1-u} L_0(a)})^{-1}.
$$
Hence for any positive $r>0$ there 
exists a constant $\kappa_2(a,r)<\infty$ such that
${\rm e}^{-sL(a,b)}\le \kappa_2(a,r)b^{-r}$. Hence
$$
(1-{\rm e}^{-sL(a,b)})^{2\ell}\ge 1-2\ell \kappa_2(a,r)b^{-r}+o(b^{-r}).
$$
This yields for any $r>0$,
\begin{eqnarray*}
\label{simplupper}
\beta (a,b) &<&  \kappa_1(t,a,r) b^{-vr}\\
&{}&  
 +\left[ 2 + \frac{1}{2} b^{-2(u-v)} + o(b^{-2(u-v)})-
 \left(1+1-b^{-(u-v)}+o(b^{-(u-v)})\right)\times \right.\\
&{}&\;\; \times \left(1-\Lambda([W'])b^{-u}+o(b^{-u})\right)
 \left(1-2\ell \kappa_2(a,r)b^{-r}+o(b^{-r})\right) \Biggr].
\end{eqnarray*}
Therefore for $\eta=\min\{vr, 2(u-v),u,r\}$
there exists $\kappa_3(t,a,\eta)$ such that $\beta (a,b)\le \kappa_3(t,a,\eta) b^{-\eta}$.
Since $r$ can be taken as big as needed, and $2(u-v)$ and $u$ can be chosen 
to be bigger that any $1-\epsilon$ for every $\epsilon\in (0,1)$
the result is shown. \hfill $\Box$

\section{Final Comments}

In \cite{hls}, (2.9) and '{Condition $\beta (\delta )$}', 
conditions on the decay of $\beta (a,b)$ are formulated which seem to be essential for the proof of limit theorems. 
The upper bound given in our Theorem \ref{decaypoly} above does not satisfy these stronger conditions.
Thus it is an open problem, whether the upper bound for STIT tessellations can be improved substantially.

\bigskip

\noindent {\bf Acknowledgments} 
The authors are indebted for the support of Program Basal CMM 
from CONICYT (Chile) and by DAAD (Germany).

\end{document}